\documentstyle{amsppt}
\def \al{\alpha}
\def \bt{\beta}
\def \gm{\gamma}

\def \KK{\Bbb {K}}
\def \QQ{\Bbb {Q}}
\def \NN{\Bbb {N}}
\def \ZZ{\Bbb {Z}}

\def \ab{\frak{ab}}

\def \sd #1#2{{#1}^{\scriptstyle{(#2)}}}

\def \as#1{\Cal P_{#1}}

\def \mm{\Cal M_S^+}

\def \asm#1{\Cal P^{\scriptscriptstyle{ab}}_{#1}}

\def\hom{\mathop{\text{\rm Hom}}}

\def \pd #1#2{#1^{(#2)}}

\def \gna #1#2{\Gamma_{#1}(#2)^{\scriptstyle {ab}}}

\def\ns#1{\vskip10pt\np{\bf #1}\qquad}
\hoffset-8pt\voffset20pt

\font\fintr=cmr9

\font\smc=cmtt10

\magnification=\magstep1

\hsize=6.5truein

\vsize=8.5truein

\document

\baselineskip=18pt

 \NoBlackBoxes

 \def\np{\par\noindent}

\def\endemo{\hfill\qed\enddemo}

\topmatter
\title On the Invariants of Matrices and the Embedding Problem
\endtitle
\author F. Vaccarino\endauthor
\address Dipartimento di  Matematica - Politecnico di Torino
- Corso Duca degli Abruzzi 24 - 10129 - Torino - Italy \quad
e-mail:vaccarino\@ syzygie.it
\endaddress
\endtopmatter
\vskip20pt

{\fintr \smallskip\leftskip=25pt \rightskip=25pt \baselineskip8pt 

\centerline{\smc{Abstract}}
Let $\KK$ be an infinite field and let $R$ be a $\KK$-algebra endowed with a homogeneous polynomial norm $N$ of 
degree $n$. We will show that $R$ is a quotient of a ring of invariants 
of $n\times n$ matrices, if $N$ satisfies a formal analogue
of the Cayley-Hamilton Theorem. 
To achieve this we use recent results on invariants of matrices in 
positive characteristic, due to S.\ Donkin and A.\ Zubkov. This is a 
first step towards proving a characteristic-free generalization of 
the main theorem proved in {\cite{8}}, which gives sufficient conditions for an 
algebra over a characteristic zero field to be embeddable into a ring 
of matrices of order $n$ over a commutative ring.
\par}
\vskip20pt
\head Introduction \endhead
\smallskip
Let $\KK$ be an infinite field and let $R$ be a $\KK$-algebra, we
would like to find sufficient conditions for the existence of a commutative $\KK$-algebra, say $L$,
such that $R$ is embeddable into $M_n(L)$, the ring of $n\times
n$ matrices over $L$. 

It is known (see {\cite{7}, \cite{8}, \cite{9}}) that, given any
ring $R$ there exist a unique commutative ring $A_n(R)$ and a unique ring homomorphism 
$j_{n}^{\scriptstyle{R}}:R\longrightarrow M_n(A_n(R))$ such that, for
any commutative ring $B$, $\hom(R,M_n(B))$ is in
bijection with $\hom(A_n(R),B)$ via $j_{n}^{\scriptstyle{R}}$. We 
will recall the construction of $(A_{n}(R), j_{n}^{\scriptstyle{R}})$ 
in the third section.

It is clear that $R$  can be embedded into a ring of  
$n\times n$ matrices over a commutative ring if and only if $j_{n}^{\scriptstyle{R}}$ 
is a monomorphism. 

The following result is proved in {\cite {8}} .
\proclaim{Theorem (C.\ Procesi)} If $R$ is a $\QQ$-algebra with trace satisfying the 
$n$-th Cayley-Hamilton identity and $j_n^{\scriptstyle{R}}:R\longrightarrow M_n(A_n(R))$ 
is the universal map we have that
$$j_n^{\scriptstyle{R}}:R\longrightarrow M_n(A_n(R))^{\scriptscriptstyle{GL_{n}(\QQ)}}$$
is an isomorphism.\endproclaim

The proof goes as follows: first, any $\QQ$-algebra $R$, with trace, satisfying the $n$-th
Cayley-Hamilton identity is showed to be of the form $D/J$, where $D$
is a ring 
of polynomial 
covariants of a suitable number of copies of the ring of
$n\times n$ matrices over $\QQ$ and $J$ is closed under trace. Then the 
surjectivity of $j_{n}^{\scriptstyle{R}}$
follows  easily from the linear reductivity of the general linear group over a 
characteristic zero field. The injectivity of the universal map is
then proved by 
means of
the Reynolds operator and using the stability of $J$ under the trace.  
\smallskip

Let $R$ be a ring endowed with a homogeneous polynomial norm $N$ of 
degree $n$. Let $t$ be a variable commuting with the elements of $R$.
Since $N$ is polynomial we can define
$$\chi_r(t):=N(t-r),$$
for all $r\in R$. 

If $\chi_r(r)=0$, for all $r\in R$, we say that $R$ is a {\it{Cayley-Hamilton ring of 
degree}} $n$.   

By means of Newton's formulas 
it follows that any $\QQ$-algebra with trace satisfying
the  $n$-th Cayley-Hamilton identity is a Cayley-Hamilton ring of 
degree $n$. 

C.\ Procesi asked the author the following natural question. 
\proclaim{Question (C.\ Procesi)} If $R$ is a Cayley-Hamilton $\KK$-algebra of degree $n$, is 
$$j_n^{\scriptstyle{R}}:R\longrightarrow M_n(A_n(R))^{\scriptscriptstyle{GL_{n}(\KK)}}$$
an isomorphism?\endproclaim

The main result of this paper is to show that any Cayley-Hamilton
$\KK$-algebra of degree $n$ is a quotient of the ring of polynomial 
covariants of a number of copies of $n\times n$ matrices over $\KK$. 

The paper is divided into three sections.
In the first section we develop some generalities on polynomial laws
and divided powers rings. 
In the second section we will prove the main result. In order to do
this we 
link the theory of invariants of matrices to the one of 
polynomial mappings and 
polynomial identities. We use the Donkin-Zubkov Theorem on invariants 
of matrices in positive characteristic  that gives a 
presentation of the ring of polynomial invariants of several matrices in 
positive characteristic.

In the last section a possible strategy for tackling the Question is outlined.

\vskip20pt
\head \S 0 Conventions and Notations \endhead
\vskip10pt

\noindent
Except otherwise stated all the rings (algebras over a
commutative ring) should  be understood associative and
with multiplicative identity.

\noindent
We denote by $\NN$ the set of non-negative integers and by $\ZZ$ the 
 ring of integers. 

\noindent
Let $R$ be a commutative ring and let $n$ be any positive integer: 
we  denote by $M_{n}(R)$ the ring of $n\times n$ matrices
with entries  in $R$. 

\noindent
Let $A$ be a set, we denote by ${\scriptstyle{\#}} A$ its
cardinality. 

\noindent
For $A$ a set and any additive monoid $M$, we denote by
$M^{\scriptstyle{(A)}}$ the  set of functions $f:A@>>> M$
with finite support.  

\noindent
Let $\al\in \sd M A$, we denote by $\mid \al \mid$ the
(finite)  sum $\sum_{a\in A} \al(a)$.

\noindent
Let $S$ be a set. We put 

$F_S:=\ZZ\langle x_s \rangle_{s\in S}$ for the free ring over 
$S$;

$F_{S}^{+}$ for the augmentation ideal of $F_{S}$ (non-unital 
subring);

$\mm$ for the free semigroup generated by $\{x_s\}_{s\in S}$.

\noindent
Let $a\in \mm$ then $d_s(a)$ denotes the degree of $a$ in $x_s$,
$\ell(a):=\sum_{s\in S}d_s(a)$ its length or total degree and 
$m(a):=(d_s(a))_{s\in S}$ its multidegree (once a  well-ordering is fixed on $S$);

\noindent
Let $a,b\in\mm$, we set $a \equiv b$ if and only if $a$ is obtained
from $b$ by a cyclic permutation of letters. This gives an equivalence
relation on $\mm$ and we denote by $\frak M_S^+$ the set of its 
equivalence classes.

\vskip20pt
\head \S 1 Multiplicative Polynomial Laws \endhead
\vskip10pt
\ns{1.1}
Let $\KK$ be a commutative ring. 
Let us recall the definition of a kind of map between $\KK$-modules
that generalizes the concept of polynomial map between free 
$\KK$-modules (see {\cite{10}} or {\cite{12}}).

\definition{Definition 1.1.1 }
Let $A$ and $B$ be two $\KK$-modules. A {\it{polynomial law }}
$\varphi$ from $A$ to $B$ is a family of mappings 
$\varphi_{_{L}}:L\otimes_{\KK} A \longrightarrow L\otimes_{\KK} B$, 
with $L$ varying in the family of  commutative $\KK$-algebras, such 
that the following diagram commutes:

$$
\CD L\otimes_{\KK}A @> \varphi_{_{L}} >> L\otimes_{\KK} B \\
    @V f\otimes 1_{_{A}}VV           @VV f \otimes 1_{_{B}}V \\
    M \otimes_{\KK} A @> \varphi_{_{M}}>> M\otimes_{\KK}B,
    \endCD
    $$ 
for all $L$, $M$ commutative $\KK$-algebras and all  
homomorphisms of $\KK$-algebras $f:L \longrightarrow M$.
\enddefinition
\smallskip
\definition{Definition 1.1.2}
Let $n\in \NN$, if $\varphi_{_{L}}(au)=a^n\varphi_{_{L}}(u)$, 
for all $a\in L$, $u\in L\otimes_{\KK} 
A$ and all commutative $\KK$-algebras $L$, then $\varphi$ will 
be said {\it{homogeneous of degree $n$}}.
\enddefinition
\smallskip

\definition{Definition 1.1.3}
If $A$ and $B$ are two $\KK$-algebras  and 
$$\cases \varphi_{_{L}}(xy)&=\varphi_{_{L}}(x)\varphi_{_{L}}(y)\\
         \varphi_{_{L}}(1_{_{L\otimes A}})&=1_{_{L\otimes B}},
         \endcases
         $$
for all commutative $\KK$-algebras $L$ and for all 
$x,y\in L\otimes A$, then $\varphi$ is called 
{\it{multiplicative}}.
\enddefinition
\smallskip

Let $A$ and $B$ be two $\KK$-modules and $\varphi:A@>>>B$ be a 
polynomial law. 
We recall the following result on polynomial laws, which is a 
restatement of Th\'eor\`eme I.1 of {\cite{10}}.
\proclaim{Proposition 1.1.4} Let $S$ be a set.
\roster
\item Let $L=\KK\otimes F_{S}$ and let $\{a_{s}\, :\,s\in S\}\subset A$
be such that $a_{s}=0$ except for a finite number of $s\in S$, then
there exist $\varphi_{\xi}((a_{s}))\in B$, with $\xi \in \sd {\NN} S$,
such that:
$$\varphi_{_{L}}(\sum_{s\in S} x_s\otimes a_{s})=\sum_{\xi \in
\NN^{(S)}}  x^{\xi}\otimes \varphi_{\xi}((a_{s})),$$ where 
$x^{\xi}:=\prod_{s\in S} 
x_s^{\xi_s}$.
\item Let $R$ be any commutative $\KK$-algebra and let 
$(r_s)_{s\in S}\subset R$, then:
$$\varphi_{_{R}}(\sum_{s\in S} r_s\otimes a_{s})=\sum_{\xi \in
\NN^{(S)}}  r^{\xi}\otimes \varphi_{\xi}((a_{s})),$$ where 
$r^{\xi}:=\prod_{s\in S} 
r_s^{\xi_s}$.
\item If $\varphi$ is homogeneous of degree $n$, then in the previous 
sum one has $\varphi_{\xi}((a_{s}))=0$ if $\mid \xi \mid$ is 
different from $n$. That is:
$$\varphi_{_{R}}(\sum_{a\in A} r_a\otimes a)=\sum_{\xi \in 
\NN^{(A)},\,\mid \xi \mid=n} r^{\xi}\otimes \varphi_{\xi}((a)).$$
In particular, if $\varphi$ is homogeneous of degree $0$ or $1$, then 
it is constant or linear, respectively.
\endroster
\endproclaim
\smallskip
Let $S$ be a set, Proposition 1.1.4 means that a polynomial law 
$\varphi:A@>>>B$ is 
completely determined by its coefficients $\varphi_{\xi}((a_{s}))$, with
$(a_s)_{s\in S} \in \sd A S$.  

\remark{Remark 1.1.5}
If $A$ is a free $\KK$-module and $\{a_{t}\, :\, t\in T\}$ is a basis 
of $A$, then $\varphi$ is completely determined by its coefficients 
$\varphi_{\xi}((a_{t}))$, with $\xi \in \NN^{(T)}$.
If also $B$ is a free $\KK$-module with basis $\{b_{u}\, :\, u\in 
U\}$, then $\varphi_{\xi}((a_{t}))=\sum_{u\in U}\lambda_{u}(\xi)b_{u}$. 
Let $a=\sum_{t\in T}\mu_{t}a_{t}\in A$. Since only a finite number of 
$\mu_{t}$ and $\lambda_{u}(\xi)$ are 
different from zero, the following makes sense:
$$\varphi(a)=\varphi(\sum_{t\in T}\mu_{t}a_{t})=
\sum_{\xi\in \NN^{(T)}} \mu^{\xi}\varphi_{\xi}((a_{t}))=
\sum_{\xi\in \NN^{(T)}} \mu^{\xi}(\sum_{u\in 
U}\lambda_{u}(\xi)b_{u})=\sum_{u\in U}(\sum_{\xi\in 
\NN^{(T)}}\lambda_{u}(\xi) \mu^{\xi})b_{u}.$$
Hence, if both $A$ and $B$ are free $\KK$-modules, a polynomial law 
between them is simply a polynomial map.
\endremark
\vskip10pt
\ns{1.2}
Let $\KK$ be any commutative ring with identity. 
For a $\KK$-module $M$ let $\Gamma(M)$ denote its divided
powers algebra (see {\cite{2}}, {\cite{10}} and {\cite{12}}). 
This is a unital commutative $\KK$-algebra, with generators $\pd m k $, 
with $m\in M$, $k \in \ZZ$ and relations, for all $m,n\in M$: 
$$
\align 
\pd m i &= 0, \qquad \forall i<0; \tag i \\
\pd m 0 &=1_{\scriptstyle{\Bbb K}}, \qquad \forall m\in M; \tag ii \\
\pd {(rm)} i &= r^i\pd m i, \qquad \forall r\in R, \forall i\in \NN; 
\tag iii \\
\pd {(m+n)} k &= \sum_{i+j=k}\pd m i \pd n j , \qquad \forall k\in
\NN; 
\tag iv \\
\pd m i \pd m j &= {i+j\choose i}\pd m {i+j} , 
\qquad \forall i,j\in \NN. \tag v 
\endalign
$$

The $\KK$-module $\Gamma(M)$ is generated by products 
(over arbitrary index sets $I$) $\prod_{i\in I} \pd {x_i} {\al_i}$
of the above generators, 
it is clear that $\prod_{i\in I} \pd {x_i} {\al_i}=0$ if $\al_i<0$ 
for some $i\in I$.  
The divided powers algebra $\Gamma(M)$ is a $\NN$-graded algebra
with homogeneous components $\Gamma_k:=\Gamma_k(M)$, ($k\in \NN$), the
submodule generated by 
$\{\prod_{i\in I} \pd {x_i} {\al_i}\; :\;\vert \al\vert=k\}$. Note
that $\Gamma_0\cong \KK$ and $\Gamma_1\cong M$. $\Gamma$ is a functor from
$\KK$-modules to commutative unital graded $\KK$-algebras.

Indeed for any morphism of $\KK$-modules $f:M @>>> N$ there 
exists a unique morphism
of graded  $\KK$-algebras $\Gamma(f):\Gamma(M) @>>> \Gamma(N)$ 
such that
$\Gamma(f)(\pd x n)=\pd {f(x)} n$, for any $x\in M$ and $n\geq 0$. 
From this it follows easily that $\Gamma$ is exact.

Furthermore
$\Gamma(L\otimes_{\KK} M)\cong L\otimes_{\KK} \Gamma(M)$ 
as graded rings by means
of $\pd {(1\otimes x)} n \mapsto 1\otimes \pd x n$.  
Thus the map $\Gamma(f)$ commutes with extensions of scalars. 

If $A$ is a (unital) $\KK$-algebra, 
then $\Gamma_k(A)$ is a (unital) $\KK$-algebra too (see
{\cite{11}}).  
To distinguish the new multiplication on $\Gamma_{k}(A)$ from the one of 
$\Gamma(A)$, we denote it by \lq\lq$\tau_k$\rq\rq. We have:
$$
\align
\prod_{i\in I}\pd {a_i} {\al_i}\;\tau_{k}\; 
\prod_{j\in J}\pd {b_j} {\beta_j}&:=
(\prod_{i\in I}\pd {a_i} {\al_i})\; \tau_k \; 
(\prod_{j\in J}\pd {b_j} {\beta_j})\\
&:=\sum_{(\lambda_{ij})\in M(\al,\beta)}\quad 
\prod_{(i,j)\in I\times J}\pd {(a_i b_j)}
{\lambda_{ij}}, \endalign
$$
where $M(\al,\beta):=\{(\lambda_{ij})\in \NN^{(I\times J)} : \sum_{i\in
I}\lambda_{ij}=\beta_j\;,\forall j\in J\; ;\sum_{j\in
J}\lambda_{ij}=\al_i \;,\forall i\in I\}$ and 
$\prod_{i\in I}\pd {a_i} {\al_i}$, 
$\prod_{j\in J}\pd {b_j} {\beta_j}\in \Gamma_k(A)$.
\smallskip

Let us denote by $\gamma_n:=(\gamma_{n,L})$ the polynomial law given by 
the composition $L\otimes M @>>> \Gamma_n(L\otimes M) @>>> L\otimes 
\Gamma_n(M)$, then $\gamma_n$ is homogeneous of degree $n$.

There is a property proved by Roby in {\cite{11}}, which 
motivates our introduction of divided powers.
\proclaim{Theorem 1.2.1}
Let $A$ and $B$ be two $\KK$ -algebras. 
The set of homogeneous multiplicative polynomial laws of 
degree $n$ from $A$ to $B$ is in bijection with the set of 
all homomorphisms of $\KK$-algebras from $\Gamma_{n}(A)$ to $B$. 
Namely, given any homogeneous multiplicative polynomial law $f:A @>>> 
B$ of degree $n$, 
there exists a unique homomorphism of $\KK$-algebras 
$\phi:\Gamma_n(A) @>>> B$ 
such that $f_{_{L}}=(1_{_{L}}\otimes \phi)\cdot \gamma_{n,L}$, 
for any commutative $\KK$-algebra $L$.
\endproclaim
\vskip10pt
\ns{1.3}
Let $B$ be a commutative ring and let 
$M_{n}(B)$ be the ring of $n\times n$ 
matrices over $B$. 
Let $b\in M_{n}(B)$ and denote by $e_{i}(b)$ the $i$-th 
coefficient of the characteristic polynomial of $b$, i.e. the trace (up to sign) of 
$\wedge^{i}(b)$. 

Let $R$ be a ring, we denote by $(R)^{\scriptscriptstyle{ab}}$ its 
abelianization, that is, its quotient by the ideal generated by the 
commutators of its elements. 

The following can be found in {\cite{13}}.
\proclaim{Proposition 1.3.1}
The ring $\gna n {M_{n}(B)}$ is isomorphic to $B$. The canonical 
projection $\ab_{n}:=\Gamma_{n}(M_{n}(B)) @>>> \gna n 
{M_{n}(B)}$ is such that, for all $b\in M_{n}(B)$ and $0\leq i \leq n$, 
$$\ab_{n}(\pd 1 {n-i} \pd b i)=e_{i}(b).$$
\endproclaim   

\ns{1.4}
Let $M$ be a free $\KK$-module and let $\{m_{i} \, : \, i\in I\}$ be a
$\KK$-basis of $M$, $I$ 
a set. Then $\Gamma_n(M)$ is a free $\KK$-module and 
$\{\prod_{i\in i}\pd {m_i} {\al_i}\,:\,
 \al \in \sd {\NN} I {\text{ with }}
\mid \al \mid =n \}$ is a $\KK$-basis of it. Now, $F_S$ is a free 
$\ZZ$-module with basis $\{1\}\cup \mm$, thus 
$$\Cal B_n:=\{\pd 1
{n-\vert\al\vert}
\prod_{\mu \in \mm}\pd {\mu}{\al_{\mu}} \,:\, 
\al \in \sd {\NN} {\mm} {\text{ with }}
\mid \al \mid \leq n \}$$ 
is a $\ZZ$-basis of $\Gamma_n(F_S)$ and we will refer to this as the
{\it{standard basis}}.

Another remark is that the natural multidegree of $F_S$ induce another
on $\Gamma_n(F_S)$ defining $m(\pd 1 {n-\vert\al\vert}
\prod_{\mu \in \mm}\pd {\mu}{\al_{\mu}}):=\sum_{\mu \in
\mm}\al_{\mu}m(\mu)$ (see \S 0). The product defined above makes then 
$\Gamma_n(F_S)$ into a graded ring, with respect to this multidegree, as
can be easily cheked.

\smallskip
For further readings on these topics we refer to 
{\cite{2}}, {\cite{10}}, {\cite{11}}, {\cite{12}}  and {\cite{13}}.

\vskip20pt
\head \S 2 Divided Powers and Invariants of Matrices \endhead
\vskip10pt
\ns{2.1}
Let us introduce the following notation: let $S$ be a set, we put

\noindent
 $A_S(n):=\ZZ[x_{ij}^s]_{1\leq i,j \leq n\,, \, s\in S}$, the 
symmetric algebra on the direct sum of 
${\scriptstyle{\#}} S$ copies of $M_{n}(\ZZ)$;

\noindent
$B_S(n):=M_n(A_S(n))$; 

\noindent
$G:=Gl_n(\ZZ);$

\noindent
$C_S(n):=A_S(n)^{\scriptscriptstyle{G}}$, 
the ring of invariants 
with respect to the action of 
$G$ on $A_S(n)$ induced by simultaneous conjugation on the direct sum 
of ${\scriptstyle{\#}} S$ copies of $M_{n}(\ZZ)$;

\vskip10pt
\ns{2.2}
Let $j_n:F_S @>>> B_S(n)$ be given by
$$j_n(x_s)=\zeta_s:=
\sum_{i,j=1}^n x_{i,j}^s e_{i,j},$$
 where $(e_{i,j})_{h,k}=\delta_{i,h}\delta_{j,k}$ for $i,j,h,k=1,\dots, n$. 
The $\zeta_s$ are the so-called "generic matrices of order $n$". 

Let 
$$\bt_{n}:\Gamma_{n}(B_{S}(n)) @>>> \gna n {B_S(n)}$$ 
be the canonical projection, then $A_{S}(n)\cong \gna n {B_S(n)}$ by
Prop.1.3.1.
 
We set  
$$E_{S}(n):=\bt_n(\Gamma_{n}(j_{n}(F_{S})) \hookrightarrow 
A_{S}(n),$$ 
by Prop.1.3.1 it is 
the subring of $A_{S}(n)$ generated by the $e_{i}(j_{n}(f))$, where 
$i\in \NN$ and $f\in F_{S}$. 
By the exactness properties of $\Gamma_{n}$ and 
$(-)^{\scriptscriptstyle{ab}}$ there exists a unique ring 
epimorphism 
$$\pi_n: \Gamma_n(F_S)^{\scriptscriptstyle{ab}}@>>> 
E_{S}(n)$$ 
such that the following diagram commutes:
$$
\CD
\Gamma_n(F_S) @> \Gamma_{n}(j_{n})>> \Gamma_{n}(j_{n}(F_{S}))  \\
@V\ab_nVV      @VV \beta_n V \\
\Gamma_n(F_S)^{\scriptscriptstyle{ab}}@> \pi_n >> E_{S}(n).
\endCD
$$
By Prop.1.3.1 it is clear that $\pi_n(\pd 1 {n-i}\pd f
i)=e_i(j_n(f))$, for all $i\in \NN$ and $f\in F_S$.

\smallskip
\remark{Remark 2.2.1}
By the previous discussion the polynomial law $j_{n}(F_{S}) @>>> 
E_{S}(n)$ that corresponds to $\bt_{n}$, via theorem 1.2.1, 
is the restriction to $j_{n}(F_{S})$ of the determinant.
\endremark
\smallskip

Recall that $E_S(n)$ is a subring of $C_S(n)$.

The Donkin-Zubkov Theorem on invariants of matrices can 
be stated in the following way, see {\cite{4}}, {\cite{5}}, 
{\cite{14}} and {\cite{15}}. (This Theorem was firstly proved by 
S.\;Donkin and then by A.\;Zubkov in another way).

\proclaim{Theorem (Donkin-Zubkov)}
The ring $C_{S}(n)$ of polynomial invariants of the direct sum 
of ${\scriptstyle{\#}} S$ copies of $n\times n$ matrices is equal to 
$E_{S}(n)$. 
\endproclaim

Then we have a surjection $\pi_n:\gna n {F_S}@>>>E_S(n)=C_S(n)$.

Let $\delta_n:A_S(n)@>>>A_S(n-1)$ be the natural projection
given by 
$$x_{ij}^s\mapsto \cases 0 &{\text{ if }} i=n {\text{ or }} j=n\\
                         x_{ij}^s &{\text{ otherwise}}.\\
\endcases
$$        
Notice that 
$(\delta_n)_{_{C_S(n)}}\big(C_S(n)\big)=
(\delta_n)_{_{E_S(n)}}\big(E_S(n)\big)=E_S(n-1)=C_S(n-1)$ and 
denote again by $\delta_n$ its restriction to $C_S(n)$. We denote by 
$C_S$ the inverse limit, in the category of graded rings, of the
inverse system $\big(C_S(n),\delta_n\big)$.

In {\cite{15}} the following is proved.

\proclaim{Theorem (Zubkov)}
The following sequence is exact:
$$0 @>>> \langle \{e_{n+1+k}(f)\, : \, k\in \NN {\text{ and }} f\in 
F_{S}^{+} \} \rangle @>>> C_S @>\theta_n>> C_S(n)@>>>0,$$
where $\theta_n$ is the canonical projection from $C_S$ to $C_S(n)$. 
\endproclaim

We are now able to state the main result of this section.
\proclaim {Theorem 2.2.2}
The map $\pi_n:\gna n {F_S}@>>>C_S(n)$ is an isomorphism of graded rings.
\endproclaim

The proof of this splits into some lemmas.

Let us consider $\Gamma(F_S^+)=\oplus_{k\geq 0}\Gamma_k(F_S^+)$. It is
a free $\ZZ$-module with basis 
$$\Cal B:=\{\prod_{\mu\in\mm}\pd {\mu}{\al_{\mu}}\,:
\,\al\in\sd {\NN}{\mm}\}.$$
So we can define a $\ZZ$-module epimorphism
$$
\align
\sigma_n:&\Gamma(F_S^+)@>>>\Gamma_n(F_S)\\
               &\prod_{\mu\in\mm}\pd {\mu}{\al_{\mu}}\mapsto
                \pd 1 {n-\vert\al\vert}\prod_{\mu\in\mm}
                 \pd {\mu}{\al_{\mu}}\\
\endalign
$$
with $n\in\NN$.
\proclaim{Lemma 2.2.3}
\roster
\item The kernel of $\sigma_n$ is $\oplus_{h>n}\Gamma_h(F_S^+)$. 
\item There is a unique product $\tau$ on $\Gamma(F_S^+)$,
that makes it an associative graded ring with identity
$1=\prod_{\mu\in\mm}\pd{\mu}{0}$, multidegree given by 
$$m(\prod_{\mu \in \mm}\pd {\mu}{\al_{\mu}}):=\sum_{\mu \in
\mm}\al_{\mu}m(\mu)$$ 
and such that $\sigma_n$ is a graded ring 
homomorphism for all $n\in\NN$, namely, 
for all  $\prod_{\mu\in\mm}\pd {\mu}{\al_{\mu}},\,\prod_{\nu\in\mm}\pd
{\nu}{\bt_{\nu}}\in \Cal B$:
$$
\big(\prod_{\mu\in\mm}\pd {\mu}{\al_{\mu}}\big)\;\tau\;
\big( \prod_{\nu\in\mm}\pd {\nu}{\bt_{\nu}}\big):=
\sum_{\gamma\in \al\odot\bt}\prod_{\mu\in\mm}\pd{\mu}{\gamma_{_{\mu,1}}}
\prod_{\nu\in\mm}\pd{\nu}{\gamma_{_{1,\nu}}}
\prod_{(\mu,\nu)\in(\mm)^2}\pd{\mu\nu}{\gamma_{_{\mu,\nu}}},
$$
where
$\al \odot \beta$ is the set of 
$(\gm_{_{\mu,\nu}})\in \sd {\NN}{\Cal M_S \times \Cal M_S}$ such that 
$\gamma_{_{1,1}}=0$ and

$$
\cases \sum_{\nu\in\Cal M_S}\gm_{_{\mu,\nu}}=\al_{\mu}, 
{\text{ for all }} \mu\in\mm \\
\sum_{\mu\in\Cal M_S}\gm_{_{\mu,\nu}}=
\beta_{\nu}, {\text{ for all }} \nu\in\mm.\\
\endcases
$$

\endroster
\endproclaim
\demo{Proof}

(1) Follows from 
$\sigma_n\big(\oplus_{h>n}\Gamma_h(F_S^+)\big)=0$ and the fact that 
$\sigma_n$ induces a bijection between 
$$\{\prod_{\mu\in\mm}\pd {\mu}{\al_{\mu}}\,:
\,\al\in\sd {\NN}{\mm} {\text{ with }}
\mid \al \mid \leq n \}$$ 
and (see 1.4)
$$\{\pd 1
{n-\vert\al\vert}
\prod_{\mu \in \mm}\pd {\mu}{\al_{\mu}} \,:\, 
\al \in \sd {\NN} {\mm} {\text{ with }}
\mid \al \mid \leq n \}=\Cal B_n.$$
(2) Let $\prod_{\mu\in\mm}\pd {\mu}{\al_{\mu}},\,\prod_{\nu\in\mm}\pd
{\nu}{\bt_{\nu}}\in \Cal B$ and let 
$$\tilde\al:=\big(n-\vert\al\vert,(\al_{mu})_{\mu\in\mm}\big){\text{
and }}\tilde\bt:=\big(n-\vert\bt\vert,(\bt_{\nu})_{\nu\in\mm}\big),$$ 
then (see 1.2):
$$
\pd 1 {n-\vert\al\vert}\prod_{\mu\in\mm}\pd {\mu}{\al_{\mu}}\;\tau_n\;
\pd 1 {n-\vert\bt\vert}\prod_{\nu\in\mm}\pd
{\nu}{\bt_{\nu}}=
\sum_{\gamma\in M(\tilde\al,\tilde\bt)}\,\prod_{(\mu,\nu)\in\Cal
M_S^2}\pd{(\mu\nu)}{\gm_{_{\mu,\nu}}}.\eqno(*)$$
Let $\gamma\in M(\tilde\al,\tilde\bt)$. By the definition of $\tau_n$
(see 1.2) one sees that:
$$
\cases \sum_{\nu\in\Cal M_S}\gm_{_{\mu,\nu}}=\al_{\mu}, 
{\text{ for all }} \mu\in\mm \\
\sum_{\mu\in\Cal M_S}\gm_{_{\mu,\nu}}=
\beta_{\nu}, {\text{ for all }} \nu\in\mm.\\
\endcases
\eqno(**)
$$
The result follows by writing the left hand side of $(*)$ as
$$
\sum_{\gamma\in \al\odot\bt}\pd 1 {n-\vert\gm\vert}\prod_{\mu\in\mm}
\pd{\mu}{\gamma_{_{\mu,1}}}
\prod_{\nu\in\mm}\pd{\nu}{\gamma_{_{1,\nu}}}
\prod_{(\mu,\nu)\in(\mm)^2}\pd{\mu\nu}{\gamma_{_{\mu,\nu}}}
$$
where
$\al \odot \beta$ is the set of 
$(\gm_{_{\mu,\nu}})\in \sd {\NN}{\Cal M_S \times \Cal M_S}$ with 
$\gamma_{_{1,1}}=0$ and satisfying $(**)$.
\endemo
\medskip
Let us denote by $\as S$ the ring $(\Gamma(F_S^+),\tau)$.
\smallskip
\proclaim{Lemma 2.2.4}
\roster
\item
The following is an epimorphism of graded rings for all $n\in\NN$:
$$\align
\rho_n:&\Gamma_{n}(F_{s})
@>>>\Gamma_{n-1}(F_{s})\\
&\pd 1 {n-\vert \al \vert}\prod_{i\in I}\pd {a_i} {\al_i} \mapsto
\pd 1 {n-1-\vert \al \vert}\prod_{i\in I}\pd {a_i} {\al_i}.\\
\endalign
$$
where $\pd 1 {n-\vert \al \vert}\prod_{i\in I}\pd {a_i} {\al_i}\in
\Gamma_n(F_S)$.
\item
$\as S$ is the graded inverse limit of the inverse system 
$\big(\Gamma_n(F_S),\rho_n\big)$ and the $\sigma_n$ are the canonical 
projections.
\item
The following is a commutative diagram in the category of graded rings

for all $n\in\NN$:
$$
\CD \gna n {F_S}@>(\rho_n)^{\scriptscriptstyle{ab}}>> \gna{n-1}{F_S}\\
@V\pi_nVV                      @VV\pi_{n-1}V\\
C_S(n)@>\delta_n>>C_S(n-1).\\
\endCD
$$
\item There exists an epimorphism of graded rings
$\pi:\asm S @>>>C_S$ such that 
$\pi(\pd f i)=\varprojlim e_i(j_n(f))$ for all $f\in F_S^+$ and $i\in\NN$.
\endroster
\endproclaim
\demo{Proof}

(1) It is a straightforward verification.

(2) By Lemma 2.2.3 (1) we have an decreasing filtration of ideals of
$\as S$:

$ker\sigma_n \subset ker\sigma_{n-1}$ and
$\sigma_n(ker\sigma_{n-1})=ker\rho_n$. 
Observe that 
$\{\big(\pd 1 {n-\vert \al \vert}\prod_{i\in I}\pd {a_i}
{\al_i}\big)_{n\in\NN}\,:\,\al\in\sd{\NN}{\mm}\}$ is a $\ZZ$-basis of
the limit, and notice that $\big(\pd 1 {n-\vert \al \vert}\prod_{i\in
I}\pd {a_i} {\al_i}\big)_{n\in\NN}\mapsto \prod_{i\in I}\pd {a_i}
{\al_i}$ is the required isomorphism.

(3) Let $I_n:=\ab_n\big(ker\rho_n\big)$, then
$I_n=\langle\{\ab_n(\pd f n)\,:\,f\in F_S^+\}\rangle$. 
Observe that $\pi_n\big(\ab_n(\pd f n)\big)=e_n(f)\in ker\delta_n$,
for all $f\in F_S^+$. 

(4) Follows directly from (3) applying the abelianization 
functor.\endemo
\medskip
Let $S_n$ be the symmetric group on $n$ letters. It acts on the
polynomial ring $\ZZ[x_1,x_2,\dots ,x_n]$ by permuting the $x_i$'s. This
action preserve the degree and we denote by $\Lambda_n^k$ the group of
invariants of degree $k$, then
$\Lambda_n:=\ZZ[x_1,x_2,\dots ,x_n]^{\scriptstyle{S_n}}=\oplus_{k\geq
0}\Lambda_n^k$. Let $q_n: \ZZ[x_1,x_2,\dots ,x_n]@>>> 
\ZZ[x_1,x_2,\dots ,x_{n-1}]$ be given by $x_n\mapsto 0$ and
$x_i\mapsto x_i$, for $i=1,\dots,n-1$. This map sends 
$\Lambda_n^k$ to $\Lambda_{n-1}^k$.Denote by $\Lambda^k$ the 
limit of the inverse system obtained in this
way. Define $\Lambda:=\oplus_{k\geq 0}\Lambda^k$: this is the
so-called ring of {\it{symmetric functions}}. Let $k\in \NN$ and set 
$$e_k(x_1,\dots,x_n):=\sum_{i_1<i_2<\dots <i_k}x_{i_1}x_{i_2}\cdots
x_{i_k}\,\in \Lambda_n^k$$ 
to denote the $k$-th elementary symmetric
polynomial in $x_1,x_2,\dots ,x_n$. It can be showed that the 
$e_k:=\varprojlim e_k(x_1,x_2,\dots ,x_n)$ are such that
$\Lambda=\ZZ[e_1,e_2,\dots,e_k,\dots]$. 

In $\Lambda$ we have a further
operation beside the product: the {\it{plethysm}}. Let $g,f\in
\Lambda$, we say that $h\in\Lambda$ is the plethysm of $g$ by $f$ and
we denote it by $h=g\circ f$ if $h$ is obtained by substituting the
monomials appearing in $f$ at the place of the variables in $g$. In
$\Lambda$ we have two other distinguished kind of functions beside 
the elementary symmetric: the
{\it{powers sums}} and the {\it{monomial symmetric functions}}. For any
$n\in \NN$ the $n$-th power sum is $p_n:=\sum_{i\geq 1}x_i^n$; for any
$h\in \NN$ and for any $\al=(\al_1,\dots,\al_h)\in \NN^h$, such that 
$\al_1\geq \al_2\geq \dots \geq \al_h>0$ the monomial symmetric function is 
$m_{\al}:=\sum_f \prod_{i\geq 1}^h x_{f(i)}^{\al_i}$, where sum is taken
over all the injections $f:\{1,\dots,h\}@>>>\NN$. Let $g\in \Lambda$,
then it is clear that $g\circ p_n=g(x_1^n,x_2^n,\dots,x_k^n,\dots)$.
  
\proclaim{Lemma 2.2.5}
\roster
\item The ring $\as S$ is generated by the $\pd {\mu} i$ with $i\in \NN$
and $\mu\in \mm$.
\item For all $a\in F_S^+$, and $n,i\in \NN$, $\pd {(a^n)} i$ belongs
to the subring of $\as S$ generated by the $\pd a j$, for $0\leq j
\leq ni$.
Furthermore $\pd {(a^n)} i = \rho_a(e_i\circ p_n)$, where 
$\rho_a :\Lambda \longrightarrow \as S$ is the
ring homomorphism defined by 
$e_j\mapsto \pd a j$. 
\endroster
\endproclaim
\demo{Proof}

(1) Observe that, for all $\al\in \NN^k$, $a\in (\mm)^k$ and $k\in 
\NN$
$$
\pd {a_1}{\al_1}\pd {a_2}{\al_2}\cdots\pd {a_k}{\al_k}=
 \pd {a_1} {\al_1}\tau\big(\pd {a_2}{\al_2}\cdots\pd
{a_k}{\al_k}\big)-
\sum
\pd{a_1}{\gamma_{_{1,0}}}
\prod_{j=2}^k\pd{a_j}{\gamma_{_{0,j}}}
\prod_{j=2}^k\pd{a_1a_j}{\gamma_{_{1,j}}},
$$
where the sum is taken over all 
$\gamma\in (\al_1)\odot (\al_2,\dots,\al_k)$ such that 
$\sum_{j=2}^k\gamma_{_{1,j}}>0$. Then 
$$\sum_{i,j} \gamma_{ij}=\sum_i \al_i - \sum_{j=2}^{k}\gamma_{_{1,j}} 
<\sum_i \al_i =\vert \al \vert.$$
Applying this reduction process inductively on $\vert \al \vert$, 
it is possible to express 
$\pd {a_1}{\al_1}\pd {a_2}{\al_2}\cdots\pd {a_k}{\al_k}$ 
as a polynomial in the $\pd {\mu} i$, with $i\in \NN$ and 
$\mu$ a monomial in the $a_j$. 

(2) Let $a\in F_S^+$
and let $t$ be a variable commuting with the $x_s$. Let us set
$$\varphi(1+ta):=\sum_{j\geq 0}t^j\pd a j,$$
then $\varphi$ is a polynomial law $\ZZ[t]\otimes_{\ZZ}F_S^+ @>>>
\ZZ[[t]]\otimes \as
S$ which is multiplicative as can be easily checked.

Let $\epsilon$ be a primitive $n$-th root of unity. 
Since $1-t^na^n=\prod_{k=1}^n (1-\epsilon^k t a)$, 
we have 
$$
\multline
\sum_{j\geq 0}(-1)^j t^{nj}\pd{(a^n)}j=\varphi(1-t^n a^n)=
\prod_{k=1}^n \varphi(1-\epsilon^k t a)=\\ 
=\prod_{k=1}^n
\big(\sum_{j_k\geq 0}(-1)^{j_k} t^{j_k}\pd{a}{j_k}\big)=
\sum_{h\geq 0}(-1)^h t^h 
\sum_{j_1+\dots+j_n=h}\prod_{k=1}^n \pd a {j_k} 
\epsilon^{j_1+2j_2+\dots+nj_n}\\
=\sum_{h\geq 0}(-1)^h t^h 
\sum_{\al\in P_{h,n}}\big(\prod_{k=1}^n \pd a {\al_k}\big)
\;m_{\al}(\epsilon,\dots,\epsilon^{n-1}, 1), \endmultline
$$
where $P_{h,n}$ denotes the set of partitions of $h$ in at most 
$n$ parts and
$m_{\al}(x_1,\dots,x_n)$ is
the monomial symmetric polynomial associated to $\al$.

\smallskip
The elementary symmetric polynomial 
$e_i\in \Lambda_n$ satisfy 
$$ e_i(\epsilon,\dots,\epsilon^{n-1}, 1)=
\cases 0 & 1\leq i \leq n-1 
\\ 1 &  i=n.
\endcases
$$
thus
$$ m_{\al}(\epsilon,\dots,\epsilon^{n-1}, 1)=
\cases 0 & n\nmid \vert \al \vert
\\ c_{\al}\in \ZZ & {\text{otherwise}},
\endcases
$$
where $c_{\al}$ is the coefficient of the $(\vert\al\vert/n)$-th power
of $e_n$ that appears in the expansion of 
$m_{\al}$ on $e_1,\dots,e_n$ in $\Lambda_n$. Therefore
$$
\sum_{j\geq 0}(-1)^j t^{nj}\pd{(a^n)}j=\sum_{j\geq 0}(-1)^{nj} t^{nj} 
\sum_{\al\in P_{nj,n}}\big(\prod_{k=1}^n \pd a {\al_k}\big)\;c_{\al},$$
thus
$$\pd {(a^n)} i =(-1)^{(n+1)i}\sum_{\al\in P_{ni,n}}
\big(\prod_{k=1}^n
\pd a {\al_k}\big)\;c_{\al}.$$

The assertion follows since this last equality is the image by 
$\rho_a$ of the one defining the plethysm $e_i\circ p_n$. \endemo

\smallskip
A monomial in $\mm $ is called {\it{indecomposable}} if it is not a
power of another. Let $\psi\subset \frak M_S^+$ be the set of (equivalence classes with
respect to cyclic permutations of) indecomposable monomials.  
Let $\ab:\as S @>>>\asm S$ be the canonical projection. 
\proclaim{Lemma 2.2.6}
The map $\pi:\asm S@>>>C_S$ is an isomorphism. 
\endproclaim
\demo{Proof}
It is easy to see, by
induction on $i\in\NN$ that $\ab\big(\pd {(ab)} i\big)=
\ab\big(\pd {(ba)} i\big)$, for all $a,b\in F_S$ and for all $i\in \NN$.
It is well known (see {\cite{1}}) that any $\mu\in\mm$ can be written 
as a power of an indecomposable monomial after cyclic permutations. 
Thus, by Lemma 2.2.5 one
gets that $\asm S$ is generated by $\ab(\pd {\mu} i)$, with 
$i\in\NN$ and $\mu\in\mm$.

S.\,Donkin showed that $C_S$ is the free commutative ring 
$\bigotimes_{\mu \in \psi} \Lambda_{\mu}$, where each $\Lambda_{\mu}$ 
is a copy of the ring of symmetric functions (see
{\cite{4}},{\cite{5}}). 
Then we have a surjection $\theta:C_S@>>>\asm S$ given by 
$e_i(\mu)\mapsto \ab(\pd {\mu} i)$, 
with $i\in\NN$ and $\mu\in\psi$. 
The result follows since $\pi\cdot\theta=id_{_{C_S}}$.
\endemo

We are now ready to prove Theorem 2.2.2
\demo{Proof of 2.2.2}
The kernel of $\asm S @>(\sigma_k)^{\scriptscriptstyle{ab}}>>
\gna n {F_S}$ is generated by the $\ab\big(\pd f k\big)$, with $k>n$  and $f\in F_S^+$.
By the previous Lemma and Zubkov Theorem this is exactly the 
image by  $\theta$ of the kernel of the projection $C_S@>>>C_S(n)$
and we are done.\endemo
\vskip10pt
\ns{2.3}
Let $D_S(n)\subset B_S(n)$ denote the ring of polynomial mappings from 
the direct sum of ${\scriptstyle{\#}} S$ copies of 
$M_{n}(\ZZ)$ to $M_n(\ZZ)$
that are $G$-equivariant with respect to the action of simultaneous 
conjugation on the direct sum of ${\scriptstyle{\#}} S$ copies of 
$M_{n}(\ZZ)$. Then $D_S(n)$ is the ring of covariants of
${\scriptstyle{\#}}S$ copies of $M_n(\ZZ)$. 
\proclaim{Proposition 2.3.1}
Let $\ZZ\{\zeta_s\}:=j_{n}(F_{S})$ be the subring of 
$B_S(n)$ generated by the generic matrices, then we have the following
equality
$$D_S(n)= C_S(n)\otimes \ZZ\{\zeta_s\}.$$
\endproclaim
\demo{Proof}
It is clear that $C_S(n)\otimes \ZZ\{\zeta_s\}\subset D_S(n)$. Let us
add one summand to $\sd {M_n(\ZZ)} S$ and set $T:=S\cup \{t\}$, with $t$ a
symbol.
Let $f\in D_S(n)$ and let $y$ be the generic matrix corresponding to
the projection of $\sd {M_n(\ZZ)} T$ on the new summand. 
Consider $e_1(fy)\in C_T(n)$, then one gets
$$e_1(fy)=\sum_{\mu} f_{\mu}\; e_1(\mu y)=
e_{1}\big((\sum_{\mu} f_{\mu}\mu)y\big),$$ 
with $f_{\mu}\in C_S(n)$, $\mu$ a monomial in the generic matrices
$\zeta_s$, $s\in S$. By non-degeneracy of the $e_1$ (=trace) on gets
that $f=\sum_{\mu} f_{\mu}\mu\in C_S(n)\otimes \ZZ\{\zeta_s\}$.\endemo
\medskip

Let us set 
$$\eta_{n}:=(\sigma_{n})^{\scriptscriptstyle{ab}}\otimes 
id_{F_{S}}:\asm S\otimes F_{S}@>>>\gna n {F_S}\otimes F_S.$$  
For all $f\in \asm S\otimes F_{S}^{+}$ and $n\in\NN$ let 
$$\chi_{n}(f):=f^n+\sum_{i=1}^n(-1)^i 
\pd f i  f^{n-i} \in \asm S\otimes F_{S}^{+}$$ 
be the $n$-th Cayley-Hamilton polynomial calculated in $f$.  

The following is the main result of this paper.
\proclaim{Theorem 2.3.2}
The following sequence is exact
$$0@>>>\langle\{\eta_{n}(\chi_n(f))\,:\,f\in \gna n {F_S}\otimes 
F_S^+\}\rangle @>>>  \gna n {F_S}\otimes F_S @>\pi_n\otimes j_n>> D_{S}(n)@>>>0.$$
\endproclaim
\demo{Proof}
Let $f\in \asm S\otimes F_{S}$ and let $y$ be a new variable.  
Then $f\mapsto \pd {(fy)} 1$ gives an isomorphism of $\ZZ$-modules between 
$ker\big( \asm S\otimes F_{S} @>>> D_S(n)\big)$ and the 
sub-$\ZZ$-module of $ker\big(\asm {S\cup \{y\}} @>>> \gna n {F_{S\cup 
\{y\}}}\big)$ 
of the elements of degree $1$ in $y$. It is then 
enough to describe the latter. 

The ideal $ker\big(\asm {S\cup \{y\}} @>>> \gna n {F_{S\cup 
\{y\}}}\big)$ it is generated as a $\ZZ$-module by elements of the form $g \pd h {n+k}$, 
with $g\in \asm {S\cup \{y\}} $, $h\in \asm {S\cup \{y\}} \otimes
F_{S\cup\{y\}}^{+}$ and $k\geq 1$. 
We claim that its elements of degree $1$ in $y$ are of the form $\pd
{(fy)} 1$, with $f\in\langle \{\chi_{n+k}(a)\,:\, k\in \NN {\text{ and }} a\in \asm
S\otimes F_{S}^+ \}\rangle.$ 

Let $d_y(g \pd h {n+k})=1$. There are two cases: 
\noindent
(1) If $d_{y}(g)=1$ then there exists $g^{'}\in \asm S\otimes F_{S}^{+}$ 
such that $g=\pd {(g^{'} y)} 1$. Then we have 
$$g \pd h {n+k}=\pd {(g^{'} y)} 1 \pd h {n+k}=
\pd {\big(g^{'}(h\chi_{n+k-1} (h)-\chi_{n+k}(h))y\big)} 1.$$

\indent
(2) If $d_{y}(g)=0$, then $h=\sum_{i}h_{i}$, where $d_{y}(h_{i})=i$ 
for all $i\in \NN$ and $h_{1}\neq 0$.  

From $\pd { h} {n+k}=\pd 
{(\sum_{i} h_{i})} {n+k}=\sum_{\vert\al\vert=n+k}\prod_{i}\pd {h_{i}} {\al_{j}}$ one sees 
that the summand of $\pd { h} {n+k}$ having degree $1$ in $y$ is $\pd{h_{0}} {n+k-1} \pd { h_{1}} 1 $. 
Now, for all $k\in \NN$ and $a,b\in \asm S\otimes F_{S}^{+}$ one gets 
$$\pd a k \pd b 1 = \sum_{i=1}^k (-1)^k \pd {(a^{k-i}b)} 1\pd a i = \pd {(\chi_k(a)\;b)} 1.$$ 
Thus $\pd{h_{0}} {n+k-1} \pd { h_{1}}
1=\pd{(\chi_{n+k-1}(h_{0})h_{1})} 1$, and the claim is proved.
 
To prove the statement it is enough to observe that 
$ker \eta_n\subset \langle\{\chi_{n+k}(f)\,:\, k\in\NN \,
f\in \asm \otimes F_S^+\}\rangle$ and that 
$\chi_{n+k}(f)\equiv f^{n+k}\chi_n(f) (\mod F_{S}\;ker \eta_n\;F_{S})$, 
for all $h\geq 0$ and for all $f\in \asm S\otimes F_{S}^{+}$.\endemo

\vskip10pt
\ns{2.4}
\proclaim{Definition 2.4.1}
Let $R$ be a ring and let $N:R@>>>R$ be a multiplicative polynomial
law homogeneous of degree $n$, with central values. We call $N$ a
homogeneous polynomial norm. Let $t$ be a
variable commuting with the elements of $R$. We set, for all $r\in R$: 
$$\chi_r(t):=N(t-r)=t^n+\sum_{i=1}^n(-1)^iN_i(r)t^{n-i}.$$
If $\chi_r(r)=0$, for all $r\in R$, then we say that $(R,N)$ is a
{\it{Cayley-Hamilton ring of degree}} $n$ (with respect to $N$).
\endproclaim
\remark{Remark 2.4.2}
Let $A$ be a commutative ring. For any integer $n$, the ring $M_n(A)$
is a Cayley-Hamilton ring of degree $n$, with respect to the usual
determinant.

The ring $D_S(n)$ is a Cayley-Hamilton ring of degree $n$, with 
respect to the usual determinant. 

The mapping
$$a\mapsto \ab_{n}(\pd a n), {\text{ for all }} a\in \gna n {F_S}\otimes F_S$$
is a polynomial norm 
$$\gna n {F_S}\otimes F_S\@>>> \gna n {F_S}$$
homogeneous of degree $n$.
Notice that $\pi_{n}\big(\ab_{n}(\pd a n)\big)=det \big(\pi_{n}\otimes 
j_{n}(a)\big)$, for all $a\in \gna n {F_S}\otimes F_S$.
\endremark
\smallskip

\proclaim{Theorem 2.4.3}
For any set $S$, the ring $D_S(n)$ is free on $S$ in the category of 
Cayley-Hamilton rings of degree $n$.\endproclaim
\demo{Proof} 
Let $R$ be a Cayley -Hamilton ring of degree $n$ with respect to $N$.
Choosen $r_s\in R$ with $s \in S$ there exists a unique ring
homomorphism 
$f:F_S \longrightarrow R$ such that $f(x_s)=r_s$. 

Let $Z(R)$ denote the center of $R$, then $N\cdot f:F_S @>>> Z(R)$
is a multiplicative polynomial law
homogeneous of degree $n$. Then there is
a unique ring homomorphism 
$$\phi:\gna n {F_S}@>>>Z(R),$$ 
such that $\phi \cdot \ab_n (\pd a n)=N\cdot  f(a)$, for all $a\in F_S$.
The map 
$$\phi\otimes f:\gna n {F_S}\otimes F_{S}@>>>R,$$ 
is then a ring homomorphism extending $f$.
Obviously this preserves the norm and its kernel must contain
$$\langle\{\eta_{n}(\chi_n(f))\,:\,f\in \gna n {F_S}\otimes 
F_S^+\}\rangle,$$ 
since $R$ is Cayley -Hamilton of degree $n$.
\endemo

\vskip20pt
\head \S 3 The Embedding Problem \endhead
\vskip10pt
\ns{3.1}
Let $\KK$ be an infinite field and let $R$ be a $\KK$-algebra, we
would like 
to find sufficient 
conditions for the existence of a commutative $\KK$-algebra, say $L$,
such that 
$R$ is embeddable into 
$M_n(L)$. 
Let us recall the construction of $(A_{n}(R), 
j_{n}^{\scriptstyle{R}})$ defined in the Introduction.

Let 
$$0 @>>> K @>>> F_S @>\pi >> R @>>>0$$ 
be any presentation of $R$ by means of generators and relations.

Let $J$ be the two sided ideal of $B_S(n)$ generated by $j(K)$. There exists an ideal $I$ of $A_S(n)$ 
such that $I=M_n(J)$. Is then easy to show (see {\cite{9}}) that 
$A_n(R)=A_S(n)/I$ and  that 
$j_{n}^{\scriptstyle{R}}(r)=j_{n}(f)+M_n(J)$, for all $r\in R$ and 
$f\in F_{S}$ such that $\pi(f)=r$.

We showed that any Cayley-Hamilton
ring of degree $n$ can be written as a quotient (in the category of
Cayley-Hamilton rings of degree $n$) of a ring of
covariants. This result can be regarded as the 
first step towards a generalization of (the proof of) Procesi's Theorem which would yield a
positive answer to the Question. 

We would like to describe a possible strategy to achieve such a
generalization. For sake 
of simplicity we put $F_S$ for $\KK\otimes F_S$, $A_S(n)$ for 
$\KK\otimes A_S(n)$ and so on.

\proclaim{Lemma 4.1.1} Let $S'\subset S$ and $I_{S'}$ be the kernel of the natural
projection  $A_S(n) \longrightarrow A_{S'}(n)$.
Let us denote by $J_{S'}$ the kernel of the induced map 
$D_{S}(n) \longrightarrow D_{S'}(n)$. It is clear that 
$J_{S'}=M_n(I_{S'})^{\scriptscriptstyle{G}}$.
The following hold:
\roster
\item $J_{S'}$ is the smallest ideal of $D_S(n)$, such that 
$e_i(J_{S'})\subset J_{S'}$, for all
$i$ and containing all the generic matrices  $\zeta_s$ with $s\notin S'$;
\item $\tilde J_{S'}:=B_S(n)J_{S'}B_S(n)$ is equal to $M_n(I_{S'})$.
\endroster
\endproclaim
\demo{Proof}

(1) Obviously $J_{S'}$ is an ideal of $D_S(n)$ such that $e_i(J_{S'})\subset J_{S'}$,
containing $\zeta_s$ with $s\notin S'$.  Let $f\in J_{S'}$, we can suppose that $f=g+h$, with
$g$ that do not contain any of $\{x_s': s\in S'\}$, while in each term of $h$ appears at
least one between $\{x_s': s\in S'\}$, as variables or in the form $e_i(MX_j)$ for certain
$i$, with $j\in S'$ and $M$ a monomial ; since one has  
$0_n=f((M_{s'})_{s'\in S'},0,\dots,0)=
g((M_s')_{s'\in S'},(M_s'')_{s''\in S-S'})=g((M_s)_{s\in S})$, 
for all $(M_s)_{s\in S}\in \sd {M_n(\KK)} S$ it follows that $g=0$ in $D(n)$. Thus $f=h$ and the assertion
follows. 

(2) It follows from (1), recalling that $J_{S'}=M_n(I_{S'})\cap 
D_{S}(n)$ and that,
by means of the elementary matrices, one can extract the entries of the matrices
belonging to $J_{S'}$, 
thus recovering the generators $\{x_{ij}^{s}: s\notin S'\}$ of $I_{S'}$
from the matrices $\zeta_s, s\notin S'$.\endemo
\smallskip

Let now $R\cong D_{S'}(n)/J$, with $e_i(J)\subset J$, for all $i=1,\dots,n$.
If $\{r_s\}_{s\in S'}$ is a set of generators of $R$ , then there 
exists a set $S$ with $S'\subset S$ and an epimorphism 
$$\rho:D_S(n)\longrightarrow D_{S'}(n),$$ 
such that $\rho(J_{S'})=J$ (it is enough to send
$\{x_s\}_{s\in S'}$ on 
$\{r_s\}_{s\in S'}$ and $x_s$, with $s\notin S'$ on the generators of
$J$). 

The above epimorphism extends to another 
$$\hat\rho:A_S(n)\longrightarrow A_{S'}(n)$$
and thus to 
$$\tilde\rho:=M_n(\hat\rho):B_S(n)\longrightarrow B_{S'}(n).$$ 
Let $I$ be an ideal of $A_{S'}(n)$, such that 
$\tilde J:=B_{S'}(n)JB_{S'}(n)=M_n(I)$.  Then $A_n(R)=A_{S'}(n)/I$. 
Notice that $\tilde\rho(\tilde J_{S'})=\tilde J$.

Since $B_{S'}(n)\cong B_{S}(n)\bigoplus \tilde J_{S'}$ as 
$GL_{n}(\KK)$-modules and $\tilde\rho$ restricted to $B_{S}(n)$ is 
the identity map, we have that $\tilde\rho$ splits and that 
$\tilde J_{S'}\cong Ker(\tilde\rho)$ by means of $a\longmapsto 
a-\tilde\rho(a)$, where we have identified $\tilde\rho$ with $i\circ 
\tilde\rho$, with $i:B_{S}(n) \hookrightarrow  B_{S'}(n)$ the 
natural embedding.

Let $V:=\tilde J_{S'}\cap Ker(\tilde\rho)$ we have the following result.
For the definition and the main theorem on good filtrations we refer the reader to
{\cite{3}}. 
\proclaim{Proposition 4.1.2}
\roster
\item If $H^{1}(GL_{n}(\KK),V)=0$ then $R \hookrightarrow 
M_{n}(A_n(R))^{GL_{n}(\KK)}$.

\item If $H^{1}(GL_{n}(\KK),V)=H^{2}(GL_{n}(\KK),V)=0$ then  $R \cong
M_{n}(A_n(R))^{GL_{n}(\KK)}$. 

\item If there is a family $\{V_{t}\}_{t\in T}$ of $G$-submodule of
$V$, with  $T$ a direct set, such that

i) each $V_{t}$ is of countable dimension and has a Good Filtration
  
ii) $V=\underset {\rightarrow} \to {\lim} \, V_{t}$,

\noindent
then $R \cong M_{n}(A_n(R))^{GL_{n}(\KK)}$.\endroster
\endproclaim
\demo{Proof}

(1) Use the short exact sequence
$$0 @>>> V @>>> \tilde J_{S'} @>>> \tilde J @>>> 0$$ 
and observe that if $H^{1}(GL_{n}(\KK),V)=0$ then 
$$0 @>>> V^{GL_{n}(\KK)} @>>> J_{S'} @>>> \tilde J^{GL_{n}(\KK)} @>>>0$$ 
is exact. But $J_{S'}/V^{GL_{n}(\KK)} \cong J\subset \tilde 
J^{GL_{n}(\KK)} $ and the statement is proved.

(2) Since 
$0=H^{1}(GL_{n}(\KK),\tilde J_{S'})\cong H^{1}(GL_{n}(\KK),Ker(\tilde\rho))$ the 
assertion follows from a) using the long cohomology sequence.

(3) If all the $V_{t}$ s have Good Filtration, then their first cohomology 
groups are zero. Since $H^{i}$ commutes with direct limits, it follows 
that $H^{1}(GL_{n}(\KK),V)=0$. Thus we have the embedding. 

Under hypotheses that are certainly 
satisfied here, the quotient of two $G$-modules having Good Filtrations has Good 
Filtration too. Hence we can apply to $\tilde J$ the observation just 
done and say that $H^{1}(GL_{n}(\KK),\tilde J)=0$. Then the 
surjectivity follows.\endemo

\vskip20pt

\centerline{\smc{Acknowledgment}}
I would like to thank C.Procesi for suggesting to work on these topics, for his
encouragement and for valuable discussions. I would like also M.Brion for useful 
remarks on this paper.

This work was done while I was a Senior Fellow at I.N.D.A.M. \lq\lq Francesco
Severi\rq\rq in Rome - Italy, I would like to thank the President and the staff of this
Institute for their support.

\enddocument
\bye